\documentclass{amsart}

\newcommand{\bi}{\bf\itshape }

\usepackage{amssymb}
\usepackage{bezier}
\usepackage{epic,eepic}
\usepackage{epsfig} 
\usepackage{times}

\date{\today}

\newcommand{\eq}[1]{~\eqref{#1}} 

\numberwithin{equation}{section} 
\newcommand{\la}{\langle} 
\newcommand{\ra}{\rangle} 

\renewcommand{\(}{\left(}
\renewcommand{\)}{\right)}

\newcommand{\R}{{\mathbb R}} 
\renewcommand{\phi}{\varphi} 
\renewcommand{\emptyset}{\varnothing} 
\newcommand{\cn}{\colon} 

\catcode`@=11 \@mparswitchfalse 

\newcounter{mnotecount}[section]

\swapnumbers

\newtheorem{thm}{Theorem}[section]
\newtheorem{lemma}[thm]{Lemma}
\newtheorem{prop}[thm]{Proposition}
\newtheorem{cor}[thm]{Corollary}

\theoremstyle{definition}

\theoremstyle{remark}

\newtheorem{remark}[thm]{Remark}

\newcommand{\f}{\partial}
\newcommand{\cd}{,\dots,} 
\renewcommand{\setminus}{\smallsetminus}
\newcommand{\s}{{\mathbb S}}

\newcommand{\id}{\operatorname{id}}
\newcommand{\lin}{\operatorname{span}}
\newcommand{\grad}{\operatorname{grad}}

\title{Smooth convex Bodies with proportional projection functions}

\author{Ralph Howard}
\address{Department of Mathematics,
University of South Carolina,
Columbia, S.C. 29208, USA}
\email{howard\char'100math.sc.edu}
\urladdr{www.math.sc.edu/$\sim$howard}

\author{Daniel Hug}
\address{Mathematisches Institut,
Albert-Ludwigs-Universit{\"a}t Freiburg,
D-79104 Freiburg,
Germany}
\email{daniel.hug@math.uni-freiburg.de}
\urladdr{http://home.mathematik.uni-freiburg.de/hug/}

\begin{document}

\begin{abstract}
For a convex body $K\subset\R^n$ and $i\in\{1,\dots,n-1\}$, the
function assigning to any $i$-dimensional subspace $L$ of $\R^n$, the
$i$-dimensional volume of the orthogonal projection of $K$ to $L$, is
called the $i$-th projection function of $K$. Let $K, K_0\subset \R^n$
be smooth convex bodies of class $C^2_+$, and let $K_0$ be centrally symmetric.
Excluding two exceptional cases, we prove that $K$ and $K_0$ are
homothetic if they have two proportional projection functions. The
special case when $K_0$ is a Euclidean ball provides an extension of
Nakajima's classical three-dimensional characterization of spheres to
higher dimensions.
\end{abstract}

\maketitle

\section{Introduction and main results}

A {\bi convex body\/} in $\R^n$ is a compact convex set with nonempty
interior. If $K$ is a convex body and $L$ a linear subspace of $\R^n$,
then $K|L$ is the orthogonal projection of $K$ onto $L$. Let
$\mathbb{G}(n,i)$ be the Grassmannian of all $i$-dimensional linear
subspaces of $\R^n$. A central question in the geometric tomography
of convex sets is to understand to what extent information about the
projections $K|L$ with $L\in \mathbb{G}(n,i)$ determines a convex
body. Possibly the most natural, but rather weak, information about
$K|L$ is its $i$-dimensional volume $V_i(K|L)$. The function
$L\mapsto V_i(K|L)$ on $\mathbb{G}(n,i)$ is the {\bi $i$-\emph{th}
projection function\/} (or the {\bi $i$-\emph{th} brightness
function\/}) of $K$. When $i=1$ this is the {\bi width function\/}
and when $i=n-1$ the {\bi brightness function\/}. If this function
is constant the body has {\bi constant $i$-brightness\/}. For
$n\geq 2$ and any $i\in
\{1,\dots,n-1\}$ by classical results about the existence of sets with
constant width and results of
Blaschke~\cite[pp.~151--154]{Blaschke:Kreis} and
Firey~\cite{Firey:constant} there are convex bodies which are not
Euclidean balls that have constant $i$-brightness
(cf.~\cite[Thm~3.3.14, p.~111; Rmk~3.3.16, p.~114]{Gardner:book}).
Thus it is not possible to determine if a convex body is a ball
from just one projection function. For other results about
determining convex bodies from a single projection function see
Chapter~3 of Gardner's book~\cite{Gardner:book} and the survey
paper~\cite{GSW97} of Goodey, Schneider, and Weil.

Therefore, as pointed out by Goodey, Schneider, and Weil
in~\cite{GSW97} and \cite{GSW97+}, it is natural to ask if a convex body with two
constant projection functions must be a ball or, more generally, what
can be said about a pair of convex bodies, one of which is centrally
symmetric, that have two of their projection functions proportional.
Examples in the smooth and the polytopal setting,
due to Campi \cite{Campi}, Gardner and Vol\v{c}i\v{c} \cite{GV93}, and to
Goodey, Schneider, and Weil \cite{GSW97+}, show that the assumption of
central symmetry on one of the bodies cannot be dropped.
Recall that a convex body is of class $C^2_+$ iff its boundary, $\f
K$, is of class $C^2$ and has everywhere positive Gauss-Kronecker curvature. A
convex body of class $C^2_+$ has a $C^2$ support function, and in fact the
class of convex bodies with $C^2$ support functions is a slightly
larger class than the class $C^2_+$. Since our proofs essentially work in
this more general class, we will consider the corresponding setting.
 A classical result~\cite{Nakajima:ball} of
S.~Nakajima (= A.~Matsumura) in 1926 states that a convex body of class
$C^2_+$ with constant width and constant brightness is a Euclidean
ball. Our main result extends this to higher dimensions:

\begin{thm}\label{Theorem1.4} Let $K, K_0\subset {\mathbb R}^{n}$ be
convex bodies with $K_0$ of class $C^2_+$ and centrally symmetric and
with $K$ having $C^2$ support function. Let $1\leq i<j\leq n-1$ be
integers such that $i\notin\{1,n-2\}$ if $j=n-1$. Assume there are
real positive constants $\alpha,\beta>0$ such that
$$
V_i(K\vert L)=\alpha V_i(K_0\vert L)\quad \text{and}\quad V_j(K\vert
U)=\beta V_j(K_0\vert U),
$$
for all $L\in\mathbb{G}(n,i)$ and $U\in
\mathbb{G}(n,j)$. Then $K$ and $K_0$ are homothetic.
\end{thm}

Other than Nakajima's result the only previously known case is $i=1$
and $j=2$ proven by Chakerian~\cite{Chakerian:rel-width} in 1967.
Letting $K_0$ be a Euclidean ball in the theorem gives:

\begin{cor}\label{intro:cor}
Let $K\subset\R^n$ be a convex body with $C^2$ support
function. Assume that $K$ has constant $i$-brightness and constant
$j$-brightness, where $1\leq i<j\leq n-1$ and $i\notin\{1,n-2\}$ if
$j=n-1$. Then $K$ is a Euclidean ball.
\end{cor}

If $\f K$ is of class $C^2$ and $K$ has constant width, then the
Gauss-Kronecker curvature of $K$ is everywhere positive. Therefore for $K$ of
class $C^2$ and of constant width the assumption of positive curvature
can be omitted:

\begin{cor}\label{Corollary1.3}
Let $K \subset {\mathbb R}^{n}$ be a convex body of class $C^{2}$ with
constant width and constant $k$-brightness for some $k \in
\{2,\dots,n-2\}$. Then K is a Euclidean ball.
\end{cor}

Unfortunately, this does not cover the case that $K$ has constant width
and brightness, which we consider the most interesting open
problem related to the subject of this paper. Under the strong
additional assumption that $K$ and $K_0$ are smooth convex bodies of revolution with
a common axis,
we can also settle the two cases not covered by Theorem \ref{Theorem1.4}.

\begin{prop}\label{revolution}
Let $K, K_0\subset\R^n$ be convex bodies that have a common axis of
revolution such that $K$ has $C^2$ support function and $K_0$ is
centrally symmetric and of class $C^2_+$. Assume that $K$ and $K_{0}$
have proportional brightness and proportional $i$-th
brightness function for an $i\in \{1,n-2\}$. Then $K$ is homothetic
to $K_0$. In particular, if $K_0$ is a Euclidean ball, then $K$ also
is a Euclidean ball.
\end{prop}

From the point of view of convexity theory the restriction to convex
bodies of class $C^2_+$ or with $C^2$ support functions is not natural
and it would be of great interest to extend Theorem~\ref{Theorem1.4}
and Corollaries~\ref{intro:cor} and~\ref{Corollary1.3} to general
convex bodies. In the case of Corollary~\ref{Corollary1.3} when
$n\geq3$, $i=1$ and $j=2$ this was done in~\cite{Howard:brightness}. 
However, from the point of view of differential geometry, the class
$C^2_+$ is quite natural and the convex bodies of constant
$i$-brightness in $C^2_+$ have some interesting differential geometric
properties. Recall that a point $x$ of $\f K$ is an {\bi umbilic
point\/} iff all of the radii of curvature of $\f K$ at $x$ are
equal. The following is a special case of
Proposition~\ref{Proposition4.4} below.

\begin{prop}\label{intro:umbilic}
Let $K$ be a convex body of class $C^2_+$ in $\R^n$ with
$n\geq 5$, and let $2\leq k\leq n-3$. Assume that $K$ has constant
$k$-brightness. Then $\f K$ has a pair of umbilic points $x_1$
and $x_2$. Moreover the tangent planes of $\f K$ at $x_1$ and $x_2$ are
parallel and the radii of curvature of $\f K$ at $x_1$ and $x_2$ are
equal.
\end{prop}

This is surprising as when $n\geq 4$ the set of $K$ in $C^2_+$ with no
umbilic points is a dense open set in $C^2_+$ with the $C^2$ topology.

Finally we comment on the relation of our results to those in the
paper~\cite{Haab:brightness} of Haab. All our main results are stated
by Haab, but his proofs are either incomplete or have errors (see the
review in Zentralblatt). In particular, the proof of his main result,
stating that a convex body of class $C^2_+$ with constant width and
constant $(n-1)$-brightness is a ball, is wrong (the proof is based
on~\cite[Lemma~5.3]{Haab:brightness} which is false even in the case
of $n=1$) and this case is still open. We have included remarks at the
appropriate places relating our results and proofs to those
in~\cite{Haab:brightness}. Despite the errors
in~\cite{Haab:brightness}, the paper still has some important
insights. In particular, while Haab's proof of his Theorem~4.1 (our
Proposition~\ref{prop:wedge}) is incomplete, see
Remark~\ref{Haab:incomplete} below, the statement is correct and is
the basis for the proofs of most of our results. Also it was Haab who
realized that having constant brightness implies the existence of
umbilic points. While his proof is incomplete and the details of the
proof here differ a good deal from those of his proposed argument, the
global structure of the proof here is still indebted to his paper.

\section{Preliminaries}\label{sec:prelim}
We will work in Euclidean space $\R^n$ with the usual inner product
$\la \cdot\,,\cdot\ra$ and the induced norm $|\cdot|$.
The support function of a convex body $K$ in $\R^n$ is the function
$h_K\cn \R^n\to \R$ given by $h_K(x)=\max_{y\in K}\la x,y\ra$. The function $h_K$
is homogeneous of degree one. A convex body is uniquely determined by
its support function. An important fact for us, first noted by
Wintner~\cite[Appendix]{Wintner:parallel}, is that if $K$ is of class $C^2_+$,
then its support function $h_K$ is of class $C^2$ on
$\R^n\setminus\{0\}$ and the principal radii of curvature (see below for a definition)
of $K$ are everywhere
positive~(cf.~\cite[p.~106]{Schneider:convex}). Conversely, if the support
function of $K$ is of class $C^2$ on $\R^n\setminus\{0\}$ 
and the principal radii of curvature of $K$ are
everywhere positive, then $K$ is of class $C^2_+$~(cf.~\cite[p.~111]{Schneider:convex}). 
In this paper, we say that a support function is of class $C^2$ if it is of class $C^2$ 
on $\R^n\setminus\{0\}$. 
Let $L$ be
a linear subspace of $\R^n$. Then the support function of the
projection $K|L$ is the restriction $h_{K|L}=h_K\big|_{L}$. In
particular, if $h_K$ is of class $C^2$, then $h_{K|L}$ is of class
$C^2$ in $L$. As an easy consequence we obtain that if $K$ is of
class $C^2_+$, then $K|L$ is of class $C^2_+$ in $L$.

All of our proofs work for convex bodies $K\subset\R^n$ that have a $C^2$ support
function, which leads to a somewhat larger class 
than the convex bodies of class  $C^2_+$. As an example, 
let $K$ be of class $C^2_+$ and let $r_0$ be the minimum of all the
radii of curvature on $\f K$. Then by Blaschke's rolling theorem
(cf.~\cite[Thm~3.2.9 p.~149]{Schneider:convex}) there is a convex set
$K_1$ and a ball $B_{r_0}$ of radius $r_0$ such that $K$ is the
Minkowski sum $K=K_1+B_{r_0}$ and no ball of radius greater than $r_0$
is a Minkowski summand of $K$. Thus no ball is a summand of $K_1$,
for if $K_1=K_2+B_r$, $r>0$, then $K=K_1+B_{r_0}=K_2+B_{r+r_0}$, contradicting
the maximality of $r_0$. As every convex body with $C^2$ boundary has
a ball as a summand, it follows that $K_1$ does not have a $C^2$
boundary. But the support function of $K_1$ is $h_{K_1}=h_K-r_0$ and
therefore $h_{K_1}$ is $C^2$. When $K_1$ has nonempty interior,
for example when $K$ is an ellipsoid with all axes of different
lengths, then $K_1$ is an example of a convex set with $C^2$ support
function, but with $\f K_1$ not of class $C^2$.

If the support function $h=h_K$ of a convex body $K\subset\R^n$ is of
class $C^2$, then let $\grad h_K$ be the usual gradient of $h_K$. This
is a $C^1$ vector field on $\R^n\setminus\{0\}$. Let $\s^{n-1}$ be the
unit sphere in $\R^n$. Then for $u\in \s^{n-1}$ the unique point on
$\f K$ with outward normal $u$ is $\grad h_K(u)$ (cf.~\cite[(2.5.8),
p.~107]{Schneider:convex}). In the case where $K$ is of class $C^2_+$,
$u\mapsto \grad h_K(u)$ is the inverse of the {\bi Gauss map\/} of $\f
K$. For this reason, $u\mapsto \grad h_K(u)$ is called the {\bi
reverse spherical image map\/}
(cf.~\cite[p.~107]{Schneider:convex}). Let $d^2h_K$ be the usual
Hessian of $h_K$ viewed as a field of selfadjoint linear maps on
$\R^n\setminus \{0\}$. That is for a vector $X$ we have that
$d^2h_KX=\nabla_X\grad h_K$ is the directional derivative of $\grad
h_K$ in the direction $X$. As $h_K$ is homogeneous of degree one for
any $u\in \s^{n-1}$, it follows that $d^2h_K(u)u=0$. Moreover, since
$d^2h_K$ is selfadjoint this implies that the orthogonal complement
$u^\bot$ of $u$ is invariant under $d^2h_K(u)$. As
$u^\bot=T_u\s^{n-1}$ we can then define a field of selfadjoint linear
maps $L(h_K)$ on the tangent spaces to $\s^{n-1}$ by
$$
L(h_K)(u):=d^2h_K(u)\big|_{u^\bot}.
$$
For given $u\in\s^{n-1}$, $L(h_K)(u)$ is called the {\bi reverse Weingarten map\/} of $K$
at $u$. The eigenvalues of $L(h_K)(u)$
are the (principal) {\bi radii of curvature\/} of $K$ in direction $u$
(cf.~\cite[p.~108]{Schneider:convex}).
Recall that if $K$ is of class $C^2_+$, then the derivative of the Gauss map
is the {\bi Weingarten map\/}
of $\f K$. As $d^2h_K$ is the directional derivative of $\grad
h_K\big|_{\s^{n-1}}$ and $\grad h_K\big|_{\s^{n-1}}$ is the inverse of
the Gauss map, we have that $L(h_K)$ is the inverse of the Weingarten
map. Provided that
$K$ is of class $C^2_+$, the Weingarten map is positive definite and therefore
the same is true of its inverse $L(h_K)$.

In the following,
the notion of the area measure of a convex body will be useful. In
the case of general convex bodies the definition is a bit involved,
see~\cite[pp.~200--203]{Schneider:convex} or
\cite[pp.~351--353]{Gardner:book}, but we will only need the case of
bodies with support functions of class $C^2$
where an easier definition is possible. Let $K\subset\R^n$ be a convex body
with support function $h_K$ of class $C^2$.
Then the (top order) {\bi area measure\/} is defined on Borel subsets
$\omega$ of $\s^{n-1}$ by
\begin{equation}\label{top-area}
 S_{n-1}(K,\omega):=\int_{\omega} \det(L(h_K)(u))\,du,
\end{equation}
where $du$ denotes integration with respect to spherical Lebes\-gue mea\-sure.
(See, for instance, \cite[(4.2.20), p.~206; Chap.~ 5]{Schneider:convex} or \cite[(A.7),
p.~353]{Gardner:book}.)

We need also a generalization of the operator $L(h_K)$. Let $K_0\subset\R^n$ be a convex
body of class $C^2_+$,
and let $h_0$ be the support function of $K_0$. As $K_0$ is of class $C^2_+$, the
linear map $L(h_0)(u)$ is positive definite for all $u\in \s^{n-1}$.
Therefore $L(h_0)(u)$ will have a unique positive definite square root
which we denote by $L(h_0)^{1/2}(u)$. Then for any convex body $K\subset\R^n$
with support function $h_K$ of class $C^2$, we define
\begin{equation}\label{def:Lh0}
L_{h_0}(h_K)(u):=L(h_0)^{-1/2}(u)L(h_K)(u)L(h_0)^{-1/2}(u)
\end{equation}
where $L(h_0)^{-1/2}(u)$ is the inverse of $L(h_0)^{1/2}(u)$. It is
easily checked that if $K$ is of class $C^2_+$, then $L_{h_0}(h_K)(u)$ is positive
definite for all $u$. Furthermore, we always have
$$
\det(L_{h_0}(h_K)(u))=\frac{\det(L(h_K)(u))}{\det(L(h_0)(u))}.
$$
The linear map $L_{h_0}(h_K)(u)$ has the interpretation as the inverse
Weingarten map in the relative geometry defined by $K_0$. This
interpretation will not be used in the present paper, but it did motivate
some of the calculations.

\section{Projections and support functions}

\subsection{Some multilinear algebra}

The geometric condition of proportional projection functions can be
translated into a condition involving reverse Weingarten maps. In
order to fully exploit this information, the following lemmas will be
used. In fact, these lemmas fill a gap in \cite[{\S}4]{Haab:brightness}.
For basic results concerning the Grassmann algebra and alternating maps, which
are used subsequently, we refer to \cite{MarcusI:73}, \cite{MarcusII:75}.

\begin{lemma}\label{Lemma4.1} Let $G,H,L\cn {\mathbb R}^{n} \to {\mathbb
R}^{n}$ be positive semidefinite linear maps. Let $k \in
\{1,\dots,n\}$, and assume that
\begin{equation}\label{eq4.1}
\left\langle \left(\land^{k}G+\land^kH\right)\xi,\xi \right\rangle =
\left\langle \left(\land^k L\right)\xi,\xi \right\rangle
\end{equation}
for all decomposable $\xi \in \bigwedge^{k}{\mathbb R}^{n}$. Then
\begin{equation}\label{eq4.1b}
\land^{k}G+\land^k H = \land^{k}L.
\end{equation}
\end{lemma}

\begin{proof} It is sufficient to consider the cases
$k \in \{2,\dots,n-1\}$. For $\xi,\zeta\in\bigwedge^k\R^n$, we define
$$
\omega_L(\xi,\zeta):=\left\langle \left(\land^k L\right)\xi,\zeta\right\rangle.
$$
Then, for any $u_{1},\dots,u_{k+1},v_{1},\dots,v_{k-1} \in {\mathbb
R}^{n}$, the identity
\begin{equation}\label{eq4.2}
\sum_{j=1}^{k+1}(-1)^{j} \omega_L(u_{1}\land\dots\land\check{u}_{j}\land\dots\land u_{k+1};
u_{j}\land v_{1}\land \dots\land v_{k-1}) = 0
\end{equation}
is satisfied, where $\check{u}_{j}$ means that $u_{j}$ is omitted. Thus, in the
terminology of
\cite{Kulkarni72}, $\omega_L$ satisfies the first Bianchi identity. Once
\eqref{eq4.2} has been verified, the proof of Lemma \ref{Lemma4.1} can be
completed as follows. Define $\omega_G$ and $\omega_H$ by replacing $L$ in the
definition of $\omega_L$ by $G$ and $H$, respectively. Then $\omega_{G,H}:=
\omega_G+\omega_H$ also satisfies the first Bianchi identity. By assumption,
$$
\omega_{G,H}(\xi,\xi)=\omega_L(\xi,\xi)
$$
for all decomposable $\xi \in \bigwedge^{k}{\mathbb
R}^{n}$. Proposition 2.1 in \cite{Kulkarni72} now implies that
$$
\omega_{G,H}(\xi,\zeta)=\omega_L(\xi,\zeta)
$$
for all decomposable $\xi,\zeta \in \bigwedge^{k}{\mathbb R}^{n}$,
which yields the assertion of the lemma.

For the proof of \eqref{eq4.2} we proceed as follows. Since $L$ is
positive semidefinite, there is a positive semidefinite linear
map $\varphi\cn {\mathbb R}^{n} \to {\mathbb R}^{n}$ such that $L =
\varphi \circ \varphi$. Hence
$$ \omega_L(u_{1}\land \dots \land u_{k}; v_{1} \land \dots \land v_{k}) = \langle \varphi
u_{1} \land \dots \land \varphi u_{k},\varphi v_{1} \land \dots
\land \varphi v_{k} \rangle $$
for all $u_{1},\dots,v_{k} \in {\mathbb R}^{n}$. For
$a_{1},\dots,a_{k+1},b_{1},\dots,b_{k-1} \in {\mathbb R}^{n}$ we
define
\begin{multline*}
\Phi(a_{1},\dots,a_{k+1}; b_{1},\dots,b_{k-1}) \\
 := \sum_{j=1}^{k+1}(-1)^{j} \langle a_{1} \land \dots \land
\check{a}_{j} \land \dots \land a_{k+1}; a_{j} \land b_{1} \land
\dots \land b_{k-1} \rangle.
\end{multline*}
We will show that $\Phi = 0$. Then, substituting $a_{i} = \varphi(u_{i})$
and $b_{j} = \varphi(v_{j})$, we obtain the required assertion \eqref{eq4.2}.

For the proof of $\Phi = 0$, it
is sufficient to show that $\Phi$
vanishes on the vectors of an orthonormal basis $e_{1},\dots,e_{n}$ of $\R^n$, since
 $\Phi$ is a multilinear map.
So let $a_{1},\dots,a_{k+1} \in \{e_{1},\dots,e_{n}\}$, whereas
$b_1,\dots,b_{k-1}$ are arbitrary.

If $a_{1},\dots,a_{k+1} $ are mutually different, then all summands of
$\Phi$ vanish, since $\langle a_{i},a_{j} \rangle = 0$ for $i \not=
j$. Here we use that
$$
\langle u_1\land\dots\land u_k,v_1\land\dots\land v_k\rangle=
\det\left(\langle u_i,v_j\rangle_{i,j=1}^k\right)
$$
for $u_1,\dots,u_k,v_1,\dots,v_k\in\R^n$.

Otherwise, $a_{i} = a_{j}$ for some $i \not= j$. In this case, we argue as follows.
Assume that $i <
j$ (say). Then, repeatedly using that $a_{i} = a_{j}$, we get
\begin{align*}
& \Phi(a_{1},\dots,a_{k+1}; b_{1},\dots,b_{k-1}) \\
& = (-1)^{i} \langle a_{1} \land \dots \land \check{a}_{i}
\land \dots \land a_{j} \land \dots \land a_{k+1}; a_{i} \land
b_{1} \land \dots \land b_{k-1} \rangle \\
 &\quad + (-1)^{j} \langle a_{1} \land \dots \land a_{i} \land
\dots \land \check{a}_{j} \land \dots \land a_{k+1}; a_{j} \land
b_{1} \land \dots \land b_{k-1} \rangle \\
 &= (-1)^{i}(-1)^{j-i-1} \langle a_{1} \land \dots \land a_{j}
\land \dots \land \check{a}_{j} \land \dots \land a_{k+1}; a_{i}
\land b_{1} \land \dots \land b_{k-1}\rangle \\
 &\quad + (-1)^{j} \langle a_{1} \land \dots \land a_{i} \land \dots
\land \check{a}_{j} \land \dots \land a_{k+1}; a_{j} \land b_{1}
\land \dots \land b_{k-1} \rangle \\
 &= 0,
\end{align*}
which completes the proof.
\end{proof}

\begin{remark}\label{Haab:incomplete}
In the proof of Theorem 4.1 in \cite{Haab:brightness}, Haab uses a
special case of Lemma \ref{Lemma4.1}, but his proof is incomplete. To
describe the situation more carefully, let $T\cn
\bigwedge^k\R^n\to\bigwedge^k\R^n$ denote a symmetric linear map
satisfying $\langle T\xi,\xi\rangle =1$ for all decomposable unit
vectors $\xi\in \bigwedge^k\R^n$. From this hypothesis Haab apparently
concludes that $T$ is the identity map (cf.\ \cite[p.\ 126, l.\
15-20]{Haab:brightness}). While Lemma \ref{Lemma4.1} implies that a
corresponding fact is indeed true for maps $T$ of a special form, a
counterexample for the general assertion is provided in \cite[p.\
124-5]{MarcusII:75}. For a different counterexample, let $k$ be even
and let $Q$ be the symmetric bilinear form defined on
$\bigwedge^k(\R^{2k})$ by $Q(w,w)=w \land w$. This is a symmetric
bilinear form as $k$ is even and $w\land w\in \bigwedge^{2k}\R^{2k}$
so that $\bigwedge^{2k}\R^{2k}$ is one dimensional and thus can be
identified with the real numbers. In this example, $Q(\xi,\xi)=0$ for
all decomposable $k$-vectors $\xi$, but $Q$ is not the zero bilinear
form.
\end{remark}

\begin{remark} Haab states a version of the next lemma,
\cite[Cor~4.2, p.~126]{Haab:brightness}, without proof.
\end{remark}

\begin{lemma}\label{Lemma4.2}
Let $G, H\cn {\mathbb R}^{n} \to {\mathbb R}^{n}$
be selfadjoint linear maps and assume that
$$
\land^{k}G+\land^kH=\beta \land^k \id
$$
for some constant $\beta\in\R$ with $\beta\neq0$ and some $k\in
\{1,\dots n-1\}$. Then $G$ and $H$ have a common orthonormal basis of
eigenvectors. If $k\geq 2$, then either $G$ or $H$ is an isomorphism.
\end{lemma}

\begin{proof} If $k=1$, this is elementary so we assume that $2\leq k\leq
n-1$. We first show that at least one of $G$ or $H$ is nonsingular.
Assume that this is not the case. Then both the kernels $\ker G$ and
$\ker H$ have positive dimension. Choose $k$ linearly independent
vectors $v_1,\dots,v_k$ as follows: If $\ker G\cap \ker H\neq \{0\}$,
then let $0\neq v_1\in \ker G\cap \ker H$ and choose any vectors
$v_2,\dots,v_k$ so that $v_1,v_2,\dots,v_k$ are linearly independent.
If $\ker G\cap \ker H= \{0\}$, then there are nonzero $v_1\in \ker G$
and $v_2\in \ker H$. Then $\ker G\cap \ker H= \{0\}$ implies that $v_1$
and $v_2$ are linearly independent. So in this case choose
$v_3,\dots,v_k$ so that $v_1,\dots,v_k$ are linearly independent. In
either case
\begin{align*}
(\land^{k}G+\land^kH)v_1\land v_2\land \dots \land v_k&
=Gv_1\land Gv_2\land \dots\land G v_k+Hv_1\land Hv_2\land \dots
\land Hv_k \\
&=0
\end{align*}
which contradicts that $\land^{k}G+\land^kH=\beta \land^k \id$ and
$\beta\neq 0$.

Without loss of generality we assume that $H$ is nonsingular. Since
$G$ is selfadjoint, there exists an orthonormal basis $e_1,\dots,e_n$
of eigenvectors of $G$ with corresponding eigenvalues
$\alpha_1,\dots,\alpha_n\in\R$. For a decomposable vector
$\xi=v_1\land\dots\land v_k\in\land^k\R^n\setminus\{0\}$, we define
\begin{align*}
\mbox{[}\xi\mbox{]}:=&\lin\{v\in\R^n:v\land\xi=0\}\\
=&\lin\{v_1,\dots,v_k\}\in\mathbb{G}(n,k).
\end{align*}
Then, for any $1\leq i_1<\dots<i_k\leq n$, we get
\begin{align*}
H(\lin\{e_{i_1},\dots,e_{i_k}\})&=\lin\{H(e_{i_1}),\dots, H(e_{i_k})\}\\
&=[H(e_{i_1})\land\dots\land H(e_{i_k})]\\
&=[(\land^k H)e_{i_1}\land\dots\land e_{i_k}]\\
&=[\left(\beta\land^k \id-\land^k G\right)e_{i_1}\land\dots\land e_{i_k}]\\
&=[(\beta-\alpha_{i_1}\cdots\alpha_{i_k})e_{i_1}\land\dots\land e_{i_k}]\\
&=\lin\{e_{i_1},\dots,e_{i_k}\},
\end{align*}
where we used that $H$ is an isomorphism to obtain the second and the last equality.
Since $k\leq n-1$, we can conclude that
\begin{align*}
H(\lin\{e_1\})&=H\left(\bigcap_{j=2}^{k+1}\lin\{e_{1},\dots,\check{e}_j,\dots,e_{k+1}\}\right)\\
&=\bigcap_{j=2}^{k+1}H\left(\lin\{e_{1},\dots,\check{e}_j,\dots,e_{k+1}\}\right)\\
&=\bigcap_{j=2}^{k+1}\lin\{e_{1},\dots,\check{e}_j,\dots,e_{k+1}\}\\
&=\lin\{e_1\}.
\end{align*}
By symmetry, we obtain that $e_i$ is an eigenvector of $H$ for $i=1,\dots,n$.
\end{proof}

\subsection{One proportional projection function}
Subsequently,
if $K,K_0\subset\R^n$ are convex bodies with support functions of class $C^2$, we
put $h:=h_K$ and $h_0:=h_{K_0}$ to simplify our notation. The following proposition
is basic for the proofs of our main results.

\begin{prop}\label{prop:wedge} Let
$K,K_0\subset \R^n$ be convex bodies having support functions of class $C^2$, let $K_0$ be
centrally symmetric, and let $k\in\{1,\dots,n-1\}$. Assume that
$\beta>0$ is a positive constant such that
\begin{equation}\label{KP}
V_k(K \vert U)=\beta V_k(K_0 \vert U)
\end{equation}
for all $U\in\mathbb{G}(n,k)$. Then, for all $u\in
\s^{n-1}$,
\begin{equation}\label{wedge}
\land^k L(h)(u)+\land^k L(h)(-u)
=2\beta\land^k L(h_0)(u) .
\end{equation}
\end{prop}

\begin{proof}
Let $u\in \s^{n-1}$ and a decomposable unit vector $\xi \in
\bigwedge^kT_u\s^{n-1}$ be fixed. Then there exist orthonormal
vectors $e_1,\dots, e_k\in u^\perp$ such that
$\xi=e_1\land\dots\land e_k$. Put $E:=\lin\{e_1\cd
e_k,u\}\in\mathbb{G}(n,k+1)$ and $E_0:=\lin\{e_1,\dots,e_k\}
\in\mathbb{G}(n,k)$. For any $v\in E\cap \s^{n-1}$,
$$
V_k\left((K\vert E)\vert (v^\perp\cap E)\right)=\beta V_k\left((K_0\vert E)\vert
(v^\perp\cap E)\right),
$$
and therefore a special case of Theorem 2.1 in \cite{GSW} (see also
Theorem 3.3.2 in \cite{Gardner:book})
 yields that
$$
S^E_k(K\vert E,\cdot)+S_k^E((K\vert E)^*,\cdot)=2\beta S_k^E(K_0\vert
E,\cdot),
$$
where $S^E_k(M,\cdot)$ denotes the (top order) surface area measure of
a convex body $M$ in $E$, and $(K\vert E)^*$ is the reflection of $K\vert E$
through the origin. Since $h_{K\vert E}=h_K\big|_E$ is of class $C^2$ in $E$,
equation~\eqref{top-area} applied in $E$ implies that
\begin{equation}\label{substi}
\det\left(d^2h_{K\vert E}(u)\big|_{E_0}\right)+
\det\left(d^2h_{K\vert E}(-u)\big|_{E_0}\right)
=2\beta\det\left(d^2h_{K_0\vert E}(u)
\big|_{E_0}\right).
\end{equation}
Since $e_1,\dots,e_k,u$ is an orthonormal basis of $E$, we further deduce that
\begin{align*}
\det\left(d^2h_{K\vert E}(u)\big|_{E_0}\right)&=\det\left(d^2h_K(u)(e_i,e_j)_{i,j=1}^k\right)\\
&=\det\left(
\langle L(h)(u)e_i,e_j\rangle_{i,j=1}^k\right)\\
&=\left\langle\land^kL(h)(u)\xi,\xi\right\rangle,
\end{align*}
and similarly for the other determinants. Substituting these expressions into
\eqref{substi} yields that
$$
\left\langle\left( \land^k L(h)(u)+\land^k L(h)(-u)\right)\xi,\xi\right\rangle=
\left\langle 2\beta \land^k L(h_0)(u)\xi,\xi\right\rangle
$$
for all decomposable (unit) vectors $\xi\in\bigwedge^k\R^n$. Hence the
 required assertion follows from Lemma \ref{Lemma4.1}.
\end{proof}

It is useful to rewrite Proposition \ref{prop:wedge} in the notation
of~\eqref{def:Lh0}. The following corollary is implied by
Proposition~\ref{prop:wedge} and Lemma~\ref{Lemma4.2}.

\begin{cor}\label{cor:L0}
Let $K,K_0\subset \R^n$ be convex bodies with $K_0$
being centrally symmetric and of class $C^2_+$ and $K$ having $C^2$
support function. Let $k\in\{1,\dots,n-1\}$. Assume that $\beta>0$ is
a positive constant such that
$$
V_k(K \vert U)=\beta V_k(K_0 \vert U)
$$
for all $U\in\mathbb{G}(n,k)$.
Then, for all
$u\in \s^{n-1}$,
\begin{equation}\label{wedge0}
\land^k L_{h_0}(h)(u)+\land^k L_{h_0}(h)(-u)
=2\beta \land^k\id_{T_u\mathbb{S}^{n-1}}.
\end{equation}
Moreover, for $k\in\{1,\dots,n-2\}$ the
linear maps $L_{h_0}(h)(u)$ and $L_{h_0}(h)(-u)$ have a common
orthonormal basis of eigenvectors.
\end{cor}

\section{The cases $1\leq i<j\leq n-2$}

\subsection{Polynomial relations}

In the sequel, it will be convenient to use the following notation.
If $x_{1},\dots,x_{n}$ are nonnegative real numbers and $I \subset
\{1,\dots,n\}$, then we put
$$ x_{I} := \prod_{\iota \in I} x_{\iota}. $$
If $I = \emptyset$, the empty product is interpreted as $x_{\emptyset} := 1$.
The cardinality of the set $I$ is denoted by $|I|$.

\begin{lemma}\label{Lemmapra}
Let $a, b > 0$ and $2 \leq k < m \leq n-1$
with $a^{m} \not= b^{k}$. Let $x_{1},\dots,x_{n}$ and $
y_{1},\dots,y_{n}$ be positive real numbers such that
$$ x_{I} + y_{I} = 2a \quad \mbox{ and } \quad x_{J} + y_{J} = 2b $$
whenever $I, J\subset\{1,\dots,n\}$, $|I| = k$ and $|J| =
m$. Then there is a constant $c>0$ such that $x_\iota/y_\iota=c$ for $\iota=1,\dots,n$.
\end{lemma}

\begin{proof} It is easy to see that this can be reduced to the case
where $m=n-1$. Thus we assume that $m=n-1$. By
assumption,
$$
x_\iota x_I+y_\iota y_I=2a\quad\text{and}\quad x_\iota x_{I'}+y_\iota y_{I'}=2a
$$
whenever $\iota\in\{1,\dots,n\}$, $I,I'\subset\{1,\dots,n\}\setminus\{\iota\}$,
$|I|=|I'|=k-1$. Subtracting these two equations, we get
\begin{equation}\label{6n}
x_\iota(x_I-x_{I'})=y_\iota(y_{I'}-y_I).
\end{equation}

By symmetry, it is sufficient to prove that $x_1/y_1=x_2/y_2$. We
distinguish several cases.

{\bf Case 1.} There exist $I,I'\subset\{3,\dots,n\}$, $|I|=|I'|=k-1$
with $x_I\neq x_{I'}$. Then \eqref{6n} implies that
$$
\frac{x_1}{y_1}=\frac{y_{I'}-y_I}{x_I-x_{I'}}=\frac{x_2}{y_2}.
$$

{\bf Case 2.} For all $I,I'\subset\{3,\dots,n\}$ with $|I|=|I'|=k-1$,
we have $x_I=x_{I'}$.

Since $1\leq k-1\leq n-3$, we obtain $x:=x_3=\dots=x_n$. From \eqref{6n}
 we get that also $y_I=y_{I'}$ for all $I,I'\subset\{3,\dots,n\}$ with
 $|I|=|I'|=k-1$. Hence, $y:=y_3=\dots=y_n$.

{\bf Case 2.1.} $x_1=x_2$. Since
$$
x_1x^{k-1}+y_1y^{k-1}=2a,\quad x_2x^{k-1}+y_2y^{k-1}=2a
$$
and $x_1=x_2$, it follows that $y_1=y_2$. In particular, we have
$x_1/y_1=x_2/y_2$.

{\bf Case 2.2.} $x_1\neq x_2$.

{\bf Case 2.2.1.} $x_1,x_2,x_3$ are mutually distinct. Choose
$$
I:=\{2\}\cup\{5,6,\dots,k+2\},\quad I':=\{4\}\cup\{5,6,\dots,k+2\}.
$$
Here note that $k+2\leq n$ and $\{5,6,\dots,k+2\}$ is the empty set for $k=2$. Then $x_I\neq x_{I'}$ as
$x_2\neq x_4=x_3$. Hence \eqref{6n} yields that
\begin{equation}\label{7n}
\frac{x_1}{y_1}=\frac{y_{I'}-y_I}{x_I-x_{I'}}=\frac{x_3}{y_3}.
\end{equation}
Next choose
$$
I:=\{1\}\cup\{5,6,\dots,k+2\},\quad I':=\{4\}\cup\{5,6,\dots,k+2\}.
$$
Then $x_I\neq x_{I'}$ as
$x_1\neq x_4=x_3$, and hence \eqref{6n} yields that
\begin{equation}\label{8n}
\frac{x_2}{y_2}=\frac{y_{I'}-y_I}{x_I-x_{I'}}=\frac{x_3}{y_3}.
\end{equation}
From \eqref{7n} and \eqref{8n}, we get $x_1/y_1=x_2/y_2$.

{\bf Case 2.2.2.} $x_1\neq x_2=x_3$ or $x_1=x_3\neq x_2$. By symmetry, it is
sufficient to consider the first case.
Since $k-1\leq n-3$ and using
$$
x_2x^{k-1}+y_2y^{k-1}=2a\quad\text{and}\quad x_3x^{k-1}+y_3y^{k-1}=2a,
$$
we get $y_2=y_3$. By the assumption of the proposition,
the equations
\begin{align}
x_2^k+y_2^k&=2a,\label{9n}\\
x_1x_2^{k-1}+y_1y_2^{k-1}&=2a,\label{10n}\\
x_2^{n-1}+y_2^{n-1}&=2b,\label{11n}\\
x_1x_2^{n-2}+y_1y_2^{n-2}&=2b. \label{12n}
\end{align}
are satisfied. From \eqref{9n} and \eqref{10n}, we get
$$
x_2^{k-1}(x_2-x_1)+y_2^{k-1}(y_2-y_1)=0.
$$
Moreover, \eqref{11n} and \eqref{12n} imply that
$$
x_2^{n-2}(x_2-x_1)+y_2^{n-2}(y_2-y_1)=0.
$$
Since $x_1\neq x_2$, we thus obtain
$$
\frac{y_1-y_2}{x_2-x_1}=\frac{x_2^{k-1}}{y_2^{k-1}}=\frac{x_2^{n-2}}{y_2^{n-2}},
$$
and therefore $y_2/x_2=1$. But now \eqref{9n}, \eqref{11n} and
$x_2=y_2$ give $x_2^k=a$ and $x_2^{n-1}=b$,
hence $a^{n-1}=b^{k}$, a contradiction. Thus Case 2.2.2 cannot occur.
\end{proof}

\begin{lemma}\label{Lemmaalge}
Let $a, b > 0$ and $1 \leq k < m \leq n-1$
with $a^{m} \not= b^{k}$. Then there exists a finite set $\mathcal{F} =
\mathcal{F}_{a,b,k,m}$, only depending on $a, b, k, m$, such that the
following is true: if $x_{1},\dots,x_{n}$ are nonnegative and
$y_{1},\dots,y_{n}$ are positive real numbers such that
$$ x_{I} + y_{I} = 2a \quad \mbox{ and } \quad x_{J} + y_{J} = 2b $$
whenever $I, J\subset\{1,\dots,n\}$, $|I| = k$ and $|J| =
m$, then $y_{1},\dots,y_{n} \in \mathcal{F}$.
\end{lemma}

\begin{remark}\label{rmk:infinite}
The condition $a^{m} \not= b^{k}$ is necessary in this lemma. For
example, if $a=b=1$, let $x_1=x_2=\dots =x_{n-1}=y_1=y_2=\dots y_{n-1}=1$,
$x_n=t$ and $y_n=1-t$, where $t\in(0,1)$. Then
$x_I+y_I=2$ for any nonemepty subset $I$ of $\{1,\dots,n\}$.
\end{remark}

\begin{proof}
It is easy to see that it is sufficient to consider the case $m=n-1$.

First, we consider the case $k=1$. Moreover, we assume that $x_1,\ldots,x_n$ are
positive. Then by assumption
\begin{equation}\label{1b}
x_\iota+y_\iota=2a\quad\text{and}\quad x_J+y_J=2b
\end{equation}
for $\iota=1,\dots,n$ and $J\subset\{1,\dots,n\}$, $|J|=n-1$. We put
$X:=x_{\{1,\dots,n\}}$ and $Y:=y_{\{1,\dots,n\}}$. Then \eqref{1b} implies
$$\frac{X}{x_\ell}+\frac{Y}{y_\ell}=2b,\quad \ell=1,\dots,n.$$
Using $y_\ell=2a-x_\ell$, this results in
$$2bx_\ell^2+(-X+Y-4ab)x_\ell+2aX=0.
$$
The quadratic equation
$$2bz^2+(-X+Y-4ab)z+2aX=0
$$
has at most two real solutions $z_1,z_2$, hence $x_1,\dots,x_n\in\{z_1,z_2\}$.

{\bf Case 1.} $x_1=\dots=x_n=:x$. Then by \eqref{1b} also $y_1=\dots=y_n=:y$. It
follows that
\begin{equation}\label{2b}
x^{n-1}+(2a-x)^{n-1}-2b=0 .
\end{equation}
The coefficient of highest degree of this polynomial equation is $2$ if $n$ is odd, and $(n-1)2a$ if
$n$ is even. Hence \eqref{2b} is not the zero polynomial.
This shows that \eqref{2b} has only finitely many solutions, which depend on $a,b,m$ only.

{\bf Case 2.} If not all of the numbers $x_1,\dots,x_n$ are equal, and hence $z_1\neq z_2$, we
 put
$$
l:=|\{\iota\in\{1,\dots,n\}:x_\iota=z_1\}|.
$$
Then $1\leq l\leq n-1$ and
 $n-l=|\{\iota\in\{1,\dots,n\}:x_\iota=z_2\}|$. Then \eqref{1b} yields that
\begin{align}
z_1^{l-1}z_2^{n-l}+(2a-z_1)^{l-1}(2a-z_2)^{n-l}&=2b,\label{3b}\\
z_1^{l}z_2^{n-l-1}+(2a-z_1)^{l}(2a-z_2)^{n-l-1}&=2b.\label{4b}
\end{align}
If $l=1$, then \eqref{3b} gives
\begin{equation}\label{5b}
z_2^{n-1}+(2a-z_2)^{n-1}=2b.
\end{equation}
Since this is not the zero polynomial, there exist only finitely many possible solutions $z_2$. Furthermore,
\eqref{4b} gives
$$z_1\left[z_2^{n-2}-(2a-z_2)^{n-2}\right]=2b-2a(2a-z_2)^{n-2}.$$
If $z_2 \neq a$, then $z_1$ is determined by this equation. The case $z_2=a$ cannot occur, since \eqref{5b}
with $z_2=a$ implies that $a^{n-1}=b$, which is excluded by assumption.

If $l=n-1$, we can argue similarly.

So let $2\leq l\leq n-2$. Note that $0<z_1,z_2<2a$ since $x_\iota,y_\iota>0$ and $x_\iota+y_\iota=2a$. Equating \eqref{3b} and
\eqref{4b}, we obtain
\begin{equation}\label{6b}
\left(\frac{2a-z_1}{z_1}\right)^{l-1}=\left(\frac{z_2}{2a-z_2}\right)^{n-l-1}.
\end{equation}
The positive points on the curve $Z_1^{l-1}=Z_2^{n-l-1}$, where $Z_1,Z_2>0$,
are parameterized by $Z_1=t^{n-l-1}$ and $Z_2=t^{l-1}$, $t>0$. Therefore setting
$$t^{n-l-1}=\frac{2a-z_1}{z_1},\qquad t^{l-1}=\frac{z_2}{2a-z_2},$$
that is
\begin{equation}\label{7b}
z_1=\frac{2a}{1+t^{n-l-1}},\qquad z_2=\frac{2at^{l-1}}{1+t^{l-1}},
\end{equation}
we obtain a parameterization of the solutions $z_1,z_2$ of \eqref{6b}. Now we substitute \eqref{7b} in \eqref{3b}
and thus get
$$
(2a)^{n-1}\frac{t^{(l-1)(n-l)}}{(1+t^{n-l-1})^{l-1}(1+t^{l-1})^{n-l}}
+(2a)^{n-1}\frac{t^{(l-1)(n-l-1)}}{(1+t^{n-l-1})^{l-1}(1+t^{l-1})^{n-l}}=2b.
$$
Multiplication by $(1+t^{n-l-1})^{l-1}(1+t^{l-1})^{n-l}$ yields a
polynomial equation where the monomial of largest degree is
$$
2b t^{(n-l-1)(l-1)}t^{(l-1)(n-l)},
$$
and therefore the equation is of degree $(l-1)(2(n-l)-1)$. This equation will have at most
$(l-1)(2(n-l)-1)$ positive solutions. Plugging these values of $t$
into \eq{7b} gives a finite set of possible solutions of \eqref{3b} and \eqref{4b},
depending only on $a,b,m$. This clearly results in a finite set of solutions of \eqref{1b}.

\bigskip

We turn to the case $2\leq k\leq n-2$. We still assume that $x_1,\ldots,x_n$ are positive.
By assumption and using Lemma \ref{Lemmapra}, we get
$$
(1+c^k)y_I=2a\quad\text{and}\quad (1+c^{n-1})y_J=2b
$$
for $I,J\subset\{1,\dots,n\}$, $|I|=k$, $|J|=n-1$, where $c>0$ is a constant such that
$x_\iota/y_\iota=c$ for $\iota=1,\dots,n$. We conclude that
$$
y_{\tilde{I}}=\frac{b}{a}\frac{1+c^k}{1+c^{n-1}}
$$
whenever $\tilde{I}\subset\{1,\dots,n\}$, $|\tilde{I}|=n-1-k$. Since $1\leq n-1-k\leq n-2$, we
obtain $y_1=\dots=y_n=:y$.
But then also $x_1=\dots=x_n=:x$.
Thus we arrive at
\begin{equation}\label{15n}
x^k+y^k=2a\quad\text{and}\quad x^{n-1}+y^{n-1}=2b.
\end{equation}
The set of positive real numbers $x,y$ satisfying \eqref{15n} is finite. In fact, \eqref{15n} implies that
$$
(2a-x^k)^{n-1}=y^{k(n-1)}=(2b-x^{n-1})^k,
$$
and thus
\begin{multline}\label{16n}
\sum_{\iota=0}^{n-1}\binom{n-1}{\iota}(2a)^\iota(-1)^{n-1-\iota}x^{k(n-1-\iota)}\\-
\sum_{\ell=0}^k\binom{k}{\ell}(2b)^\ell(-1)^{k-\ell}x^{(n-1)(k-\ell)}=0.
\end{multline}
The coefficient of the monomial of highest degree is $(-1)^{n-1}+(-1)^{k-1}$,
if this number is nonzero, and otherwise
it is equal to $(n-1)(2a)(-1)^{n-2}$, since $k(n-2)>(n-1)(k-1)$. In any case,
the left side of \eqref{16n} is not
the zero polynomial, and therefore \eqref{16n} has only a finite number of
solutions, which merely depend on $a,b,k,m$.

\medskip

Finally, we turn to the case where some of the numbers $x_1,\ldots,x_n$ are zero. For instance,
let $x_1=0$. Then we obtain that
$$
y_1y_{I'}=2a,\qquad y_1y_{J'}=2b
$$
whenever $I',J'\subset\{2,\ldots,n\}$, $|I'|=k-1$ and $|J'|=n-2$, and thus
 $y_{J'}/y_{I'}=b/a$. Therefore $y_{\tilde I}=b/a$ for all $\tilde I\subset \{2,\ldots,n\}$ with
 $|\tilde I|=n-1-k$. Using that $k\ge 1$, we find that $y:=y_2=\ldots=y_n=(b/a)^{\frac{1}{n-1-k}}$.
Since $y_1y^{k-1}=2a$, we again get that $y_1,\ldots,y_n$ can assume only finitely many values, depending
only on $a,b,k,m=n-1$.
\end{proof}

\subsection{Proof of Theorem \ref{Theorem1.4} for $1\leq i<j\leq n-2$}

An application of Corollary \ref{cor:L0} shows that, for $u\in \s^{n-1}$,
\begin{equation}\label{re1}
\land^i L_{h_0}(h)(u)+
\land^i L_{h_0}(h)(-u)=2\alpha \land^i \text{id}_{u^\perp},
\end{equation}
\begin{equation}\label{re2}
\land^j L_{h_0}(h)(u)+\land^j L_{h_0}(h)(-u)=2\beta\land^j
\text{id}_{u^\perp},
\end{equation}
Since $i<j\leq n-2$, Corollary \ref{cor:L0} also implies that, for any fixed
$u\in \s^{n -1}$, $L_{h_0}(h)(u)$ and $L_{h_0}(h)(-u)$ have a common
orthonormal basis of eigenvectors with corresponding nonnegative
eigenvalues $x_1,\dots,x_{n-1}$ and $y_1,\dots,y_{n-1}$, respectively.

{\bf Case 1.} $\alpha^j\neq\beta^i$. We will show that
there is a finite set, $\mathcal{F}^*_{\alpha,\beta,i,j}$, independent
of $u$, such that
\begin{equation}\label{re3}
\det\(L_{h_0}(h)(u)\)=
\frac{\det L(h)(u)}{\det
L(h_0)(u)}\in\mathcal{F}^*_{\alpha,\beta,i,j},\quad \text{for
all $u\in\s^{n-1}$}.
\end{equation}
Assume this is the case. Then, since $h,h_0$ are of class $C^2$, the
function on the left-hand side of \eqref{re3} is continuous on the
connected set $\s^{n-1}$ and hence must be equal to a constant
$\lambda\geq 0$. If $\lambda=0$, then $\det L(h)\equiv0$ and, as $\det
L(h)$ is the density of the surface area measure $S_{n-1}(K,\cdot)$
with respect to spherical Lebesgue measure, this implies that the surface
area measure $S_{n-1}(K,\cdot)\equiv 0$. But this cannot be true, since $K$
is a convex body (with nonempty interior). Therefore
$\lambda>0$. Again using that $\det L(h)(u)$ is the density of the
surface measure $S_{n-1}(K,\cdot)$, and similarly for $h_0$ and $K_0$,
we obtain $S_{n-1}(K,\cdot)= S_{n-1}(\lambda^{1/(n-1)}K_0,\cdot)$. But
then Minkowski's inequality implies that $K$ and $K_0$ are homothetic
(see
\cite[Thm~7.2.1]{Schneider:convex}).

To construct the set $\mathcal{F}^*_{\alpha,\beta,i,j}$, we first put
$0$ in the set. Then we only have to consider the points $u\in
\s^{n-1}$ where $\det L_{h_0}(h)(u)\neq0$. At these points \eqref{re1} and
\eqref{re2} show that the assumptions of Lemma~\ref{Lemmaalge} are
satisfied (with $n$ replaced by $n-1$). Hence there is a finite set
$\mathcal{F}_{\alpha,\beta,i,j}$, such that for any $u\in \s^{n-1}$
with $\det L_{h_0}(h)(u)\neq0$,
if $x_1,\dots,x_{n-1}$ are the eigenvalues of $L_{h_0}(h)(-u)$
and $y_1,\dots,y_{n-1}$ are the eigenvalues of $L_{h_0}(h)(u)$, then
$y_1,\dots,y_{n-1}\in \mathcal{F}_{\alpha,\beta,i,j}$.
Let $\mathcal{F}^*_{\alpha,\beta,i,j}$ be the union of $\{0\}$ with
the set of all products of $n-1$ numbers each from the set
$\mathcal{F}_{\alpha,\beta,i,j}$.

{\bf Case 2.} If $\alpha^j=\beta^i$, then the assumptions can be
rewritten in the form
\begin{equation}\label{e4.8}
\left(\frac{V_j(K_0\vert U)}{V_j(K\vert U)}\right)^{\frac{1}{j}}=
\left(\frac{V_i(K_0\vert L)}{V_i(K\vert L)}\right)^{\frac{1}{i}}
\end{equation}
for all $U\in\mathbb{G}(n,j)$ and all $L\in\mathbb{G}(n,i)$. Let
$U\in\mathbb{G}(n,j)$ be fixed. By homogeneity we can replace $K_0$ by
$\mu K_0$ on both sides of \eqref{e4.8}, where $\mu>0$ is
chosen such that $V_j(\mu K_0\vert U)=V_j(K\vert U)$. We put
$M_0:=\mu K_0\vert U$ and $M:=K\vert U$. Then, for any $L\in
\mathbb{G}(n,i)$ with $L\subset U$, we have
$$
V_j(M)=V_j(M_0)\quad\text{and}\quad V_i(M\vert L)=V_i(M_0\vert L).
$$
By the discussion in \cite[{\S} 4]{GSW97} or the main theorem in
\cite{ChLu}, we infer that $M$ is a translate of $M_0$, and therefore
$K| U$ and $K_0| U$ are homothetic. Since $j\geq 2$, Theorem 3.1.3 in
\cite{Gardner:book} shows that $K$ and $K_0$ are homothetic. \qed

\section{The cases $2\leq i<j\leq n-1$ with $i\neq n-2$}

\subsection{Existence of relative umbelics}

We need another lemma concerning polynomial relations.

\begin{lemma} Let $n\geq 5$, $k\in\{2,\dots,n-3\}$, $\gamma>0$, and let positive real numbers
$0<x_1\leq x_2\leq \dots\leq x_{n-1}$ be given.
Assume that
\begin{equation}\label{sternl}
x_I+x_{I^*}=2\gamma
\end{equation}
for all $I\subset\{1,\dots,n-1\}$, $|I|=k$, where $I^*:=\{n-i:i\in I\}$. Then $x_1=\dots=x_{n-1}$.
\end{lemma}

\begin{proof}
Choosing $I=\{1,2,\dots,k\}$ in \eqref{sternl}, we get
\begin{equation}\label{eql1}
x_1x_2\cdots x_k+x_{n-k}\cdots x_{n-2}x_{n-1}=2\gamma.
\end{equation}
Choosing $I=\{1,n-k,\dots,n-2\}$ in \eqref{sternl}, we obtain
\begin{equation}\label{eql2}
x_1x_{n-k}\cdots x_{n-2}+x_2\cdots x_kx_{n-1}=2\gamma.
\end{equation}
Subtracting \eqref{eql2} from \eqref{eql1}, we arrive at
\begin{equation}\label{eql3}
x_{n-k}\cdots x_{n-2}(x_{n-1}-x_1)+x_2\cdots x_k(x_1-x_{n-1})=0.
\end{equation}

Assume that $x_1\neq x_{n-1}$. Then \eqref{eql3} implies that
\begin{equation}\label{eql4}
x_2\cdots x_k=x_{n-k}\cdots x_{n-2}.
\end{equation}
We assert that $x_2=x_{n-2}$. To verify this, we first observe that
 $2\leq k\leq n-3$ and $x_2\leq \dots \leq x_{n-2}$.
After cancellation of factors with the same index on both sides of \eqref{eql4},
we have
\begin{equation}\label{eql5}
x_2\cdots x_l=x_{n-l}\cdots x_{n-2},
\end{equation}
where $2\leq l<n-l$ (here we use $k\leq n-3$). Since
$$
x_l\leq x_{n-l},\quad x_{l-1}\leq x_{n-l+1},\quad\dots\quad x_2\leq x_{n-2},
$$
equation \eqref{eql5} yields that $x_2=\dots=x_{n-2}$.

Now \eqref{eql1} turns into
\begin{equation}\label{eql6}
x_1x_2^{k-1}+x_2^{k-1}x_{n-1}=2\gamma.
\end{equation}
From \eqref{sternl} with $I=\{2,\dots,k+1\}$ and using that $k\leq n-3$, we obtain
\begin{equation}\label{eql7}
x_2^k+x_2^k=2\gamma.
\end{equation}
Hence \eqref{eql6} and \eqref{eql7} show that
\begin{equation}\label{eql8}
x_1+x_{n-1}=2x_2.
\end{equation}
Applying \eqref{sternl} with $I=\{1,\dots,k-1,n-1\}$ and using \eqref{eql7}, we get
$$
2x_1x_2^{k-2}x_{n-1}=2\gamma=2x_2^k,
$$
hence
\begin{equation}\label{eql9}
x_1x_{n-1}=x_2^2.
\end{equation}
But \eqref{eql8} and \eqref{eql9} give $x_1=x_{n-1}$, a contradiction.

This shows that $x_1=x_{n-1}$, which implies the assertion of the lemma.
\end{proof}

\begin{prop}\label{Proposition4.4}
Let $K,K_0 \subset \R^n$ be convex bodies with $K_0$ centrally
symmetric and of class $C^2_+$ and $K$ having a $C^2$ support
function. Let $n\geq 5$ and $k\in\{2,\dots,n-3\}$. Assume that there
is a constant $\beta>0$ such that
$$
V_k(K\vert U)=\beta V_k(K_0\vert U)
$$
for all $U\in\mathbb{G}(n,k)$. Then there exist $u_0\in \s^{n-1}$ and
$r_0>0$ such that
$$
 L_{h_0}(h)(u_0)=L_{h_0}(h)(-u_0)=r_0\id_{T_{u_0}\s^{n-1}}.
$$
\end{prop}

\begin{proof} For $u\in \s^{n-1}$, let $r_1(u),\dots,r_{n-1}(u)$ denote
the eigenvalues of the selfadjoint linear map $L_{h_0}(h)(u)\cn
T_u\s^{n-1}\to T_u\s^{n-1}$, which are ordered such that
$$
r_1(u)\leq\dots\leq r_{n-1}(u).
$$
Then we define a continuous map $R\cn
\s^{n-1}\to \R^{n-1}$ by
$$
R(u):=(r_1(u),\dots,r_{n-1}(u)).
$$
By the Borsuk-Ulam theorem (cf.\ \cite[p.~93]{Guillemin-Pollack} or
\cite{Matousek:BU}), there is some $u_0\in \s^{n-1}$ such that
\begin{equation}\label{4.0a}
R(u_0)=R(-u_0).
\end{equation}
Corollary \ref{cor:L0} shows that $L_{h_0}(h)(u_0)$ and
$L_{h_0}(h)(-u_0)$ have a common orthonormal basis
$e_1,\dots,e_{n-1}\in u_0^\perp$ of eigenvectors and that at least one
of $L_{h_0}(h)(u_0)$ or $L_{h_0}(h)(-u_0)$ is nonsingular. But
$R(u_0)=R(-u_0)$ implies that $L_{h_0}(h)(u_0)$ and $L_{h_0}(h)(-u_0)$
have the same eigenvalues and thus they are both nonsingular.
Therefore the eigenvalues of both $L_{h_0}(h)(u_0)$ and
$L_{h_0}(h)(-u_0)$ are positive.

We can assume that, for $\iota=1,\dots,n-1$, $e_\iota$ is an
eigenvector of $L_{h_0}(h)(u_0)$ corresponding to the eigenvalue
$r_\iota:=r_\iota(u_0)$. Next we show that $e_\iota$ is an eigenvector
of $L_{h_0}(h)(-u_0)$ corresponding to the eigenvalue
$r_{n-\iota}(-u_0)$. Let $\tilde{r}_\iota$ denote the eigenvalue of
$L_{h_0}(h)(-u_0)$ corresponding to the eigenvector $e_\iota$,
$\iota=1,\dots,n-1$. Since $\tilde{r}_1,
\dots,\tilde{r}_{n-1}$ is a permutation of $r_1(-u_0),\dots,r_{n-1}(-u_0)$, it is sufficient to show
that $\tilde{r}_1\geq \dots\geq \tilde{r}_{n-1}$. By Corollary
\ref{cor:L0}, for any $1\leq i_1<\dots< i_k\leq n-1$ we have
$$
\left(\land^kL_{h_0}(h)(u_0)+\land^k L_{h_0}(h)(-u_0)\right)e_{i_1}\land\dots\land e_{i_k}=2\beta
e_{i_1}\land\dots\land e_{i_k} ,
$$
and therefore
\begin{equation}\label{4.1a}
r_{i_1}\cdots r_{i_k}+\tilde{r}_{i_1}\cdots \tilde{r}_{i_k}=2\beta .
\end{equation}
For $\iota\in\{1,\dots,n-2\}$, we can choose a subset $I\subset\{1,\dots,n-1\}$ with
$|I|=k-1$ and $\iota,\iota+1\notin I$, since $k+1\leq n-1$. Then \eqref{4.1a} yields
$$
r_I r_\iota+\tilde{r}_I \tilde{r}_\iota=
r_I r_{\iota+1}+\tilde{r}_I\tilde{r}_{\iota+1}
\geq r_I r_{\iota}+ \tilde{r}_I\tilde{r}_{\iota+1},
$$
which implies that $\tilde{r}_\iota\geq\tilde{r}_{\iota+1}$.

Let $1\leq i_1<\dots<i_k\leq n-1$ and $I:=\{i_1,\dots,i_k\}$. Applying the linear map
$\land^k L_{h_0}(h)(u_0)+
\land^k L_{h_0}(h)(-u_0)$ to $e_{i_1}\land\dots\land e_{i_k}$, we get
\begin{equation}\label{eqsp1}
\prod_{\iota\in I}r_\iota(u_0)+ \prod_{\iota\in I}r_{n-\iota}(-u_0)=2\beta.
\end{equation}
From \eqref{4.0a} and \eqref{eqsp1} we conclude that the sequence
$0<r_1(u_0)\leq\dots\leq r_{n-1}(u_0)$ satisfies the hypothesis of Lemma
\ref{sternl}. Hence, $r_1(u_0)=\dots=r_{n-1}(u_0)=:r_0$. But
$R(-u_0)=R(u_0)$ implies that also
$r_1(-u_0)=\dots=r_{n-1}(-u_0)=r_0$, which yields the assertion of the
proposition.
\end{proof}

\subsection{Proof of Theorem \ref{Theorem1.4}: remaining cases}
It remains to consider the cases where $j=n-1$. Hence, we have $2\leq i\leq n-3$. Proposition
\ref{Proposition4.4} implies that there is some $u_0\in \s^{n-1}$ such that the eigenvalues of
$L_{h_0}(h)(u_0)$ and $L_{h_0}(h)(-u_0)$ are all equal to $r_0>0$. But then Corollary \ref{cor:L0}
 shows that
$$
r_0^i+r_0^i=2\alpha=2\frac{V_i(K\vert L)}{V_i(K_0\vert L)},
$$
for all $L\in\mathbb{G}(n,i)$, and
$$
r_0^j+r_0^j=2\beta=2\frac{V_j(K\vert U)}{V_j(K_0\vert U)},
$$
for all $U\in\mathbb{G}(n,j)$. Hence, we get
$$
\left(\frac{V_j(K_0\vert U)}{V_j(K\vert U)}\right)^{\frac{1}{j}}=
\left(\frac{V_i(K_0\vert L)}{V_i(K\vert L)}\right)^{\frac{1}{i}}
$$
for all $U\in\mathbb{G}(n,j)$ and all $L\in\mathbb{G}(n,i)$. Thus again
 equation \eqref{e4.8} is available and the proof can be completed as before.
\qed

\subsection{Proof of Corollary \ref{Corollary1.3}}
Let $K$ have constant width $w$. Then, \cite[{\S}64]{Bonnesen-Fenchel},
the diameter of $K$ is also $w$ and any point $x\in \f K$ is the
endpoint of a diameter of $K$. That is there is $y\in \partial K$
such that $|x-y|=w$. Then $K$ is contained in the closed ball
$B(y,w)$ of radius $w$ centered at $y$ and $x\in \f B(y,w)\cap K$.
Thus if $\f K$ is $C^2$, then $\f K$ is internally tangent to the
sphere $\f B(y,w)$ at $x$. Therefore all the principle curvatures of
$\f K$ at $x$ are greater or equal than the principle curvatures of $\f B(y,w)$ at
$x$, and thus all the principle curvatures of $\f K$ at $x$ are at
least $1/w$. Whence the Gauss-Kronecker curvature of $\f K$ at $x$ is at
least $1/w^{n-1}$. As $x$ was an arbitrary point of $\f K$ this shows
that if $\f K$ is a $C^2$ submanifold of $\R^n$ and $K$ has constant
width, then $\f K$ is of class $C^2_+$. Corollary~\ref{Corollary1.3} now
follows directly from Corollary~\ref{intro:cor}.\qed

\section{Bodies of revolution}

We now give a proof of Proposition \ref{revolution}. By assumption, there are constants
$\alpha,\beta>0$ such that
$$
V_i(K | L)=\alpha V_i(K_0 | L)\quad \text{and}\quad V_{n-1}(K |
U)=\beta V_{n-1}(K_0 | U),
$$
for all $L\in\mathbb{G}(n,i)$ and $U\in \mathbb{G}(n,n-1)$, where
$i\in\{1,n-2\}$. We can assume that the axis of revolution contains
the origin and has direction $e\in\s^{n-1}$. Let
$u\in\s^{n-1}\setminus\{\pm e\}$. Then there are $\varphi\in
\left(-\frac{\pi}{2},\frac{\pi}{2}\right)$ and $v_0\in\s^{n-1}\cap
u^\perp$ such that $u=\cos\varphi\, v_0+\sin\varphi\, e$. For the sake
of completeness we include a proof of the following lemma.

\begin{lemma}\label{eigenvalues}
The map $L(h_K)(u)$ is a multiple of the identity map on $e^\perp\cap
v_0^\perp$ and has $-\sin\varphi\, v_0+\cos \varphi\, e$ as an eigenvector.
\end{lemma}

\begin{proof} By rotational invariance, there is some
$r(\varphi)>0$ such that
\begin{equation}\label{diff}
h_K(\cos\varphi\,v+\sin\varphi\, |v|e)=r(\varphi)|v|,
\end{equation}
for all $v\in e^\perp$. Differentiating \eqref{diff} twice with
respect to $v\in e^\perp$ yields that, for any $v,w\in e^\perp\cap
v_0^\perp$,
$$
\cos^2\varphi \, d^2h_K(\cos\varphi\, v_0+\sin\varphi\, e)(v,w)
 =r(\varphi) \langle v,w\rangle.
$$
Moreover, differentiating~\eqref{diff} with respect to $v$, we obtain,
for any $v\in e^\perp\cap v_0^\perp$,
\begin{equation}\label{diffdiff}
dh_K(\cos\varphi\, v_0+\sin\varphi\, e)(v)=0.
\end{equation}
Differentiating \eqref{diffdiff} with respect to $\varphi$, we obtain
$$
d^2h_K(\cos\varphi\, v_0+\sin\varphi\, e)(v,-\sin\varphi\, v_0+\cos\varphi\, e)=0.
$$
Thus, if $v_1,\dots,v_{n-2}$ is an orthonormal basis of $e^\perp\cap
v_0^\perp$, then $-\sin\varphi\, v_1+\cos\varphi\, e,v_1,\dots,$
$v_{n-2}$ is an orthonormal basis of eigenvectors of $L(h_K)(u)$ with
corresponding eigenvalues $x_1$ and $x_2=\dots=x_{n-1}=:x$.
\end{proof}

\begin{proof}[Proof of Proposition~\ref{revolution}]
Let $K$ and $K_0$ be as in Proposition~\ref{revolution} and let $e$ be
a unit vector in the direction of the common axis of rotation of $K$
and $K_0$. Let $h$ be the support function of $K$ and $h_0$ the
support function of $K_0$. Let $u\in \s^{n-1}\cap e^\bot$ be a point
in the equator of $\s^{n-1}$ defined by $e$. As $e$ is orthogonal to
$u$, the vector $e$ is in the tangent space to $\s^{n-1}$ at $u$. Let
$e_2\cd e_{n-1}$ be an orthonormal basis for $\{u,e\}^\bot$. Then
$e,e_2\cd e_{n-1}$ is an orthonormal basis for both $T_u\s^{n-1}$ and
$T_{-u}\s^{n-1}$. By Lemma~\ref{eigenvalues} there are eigenvalues
$x_1$, and $x_2=x_3=\dots =x_{n-1}=:x$ such that $L(h)(u)e=x_1 e$ and
$L(h)(u)e_j=x e_j$ for $j=2\cd n-1$. By rotational symmetry we also
have $L(h)(-u)e=x_1 e$ and $L(h)(-u)e_j=x e_j$ for $j=2\cd n-1$.
Likewise if $y_1$, and $y_2=y_3=\dots =y_{n-1}=:y$ are the eigenvalues
of $L(h_0)(u)$, then they are also the eigenvalues of $L(h_0)(-u)$ and
$L(h_0)(\pm u)e=y_1 e$ and $L(h_0)(\pm u)e_j=y e_j$ for $j=2\cd n-1$.
Proposition~\ref{prop:wedge} implies the polynomial relations
\begin{align*}
x_1x^{i-1}+x_1x^{i-1}&=2\alpha y_1y^{i-1},\\
x^{i}+x^{i}&=2\alpha y^i,\\
x_1x^{n-2}+x_1x^{n-2}&=2\beta y_1y^{n-2}.
\end{align*}
The first two of these imply that $x/y=x_1/y_1$ and therefore
$$
\alpha^{n-1}=\(\frac{x}{y}\)^{i(n-1)}=\beta^i.
$$
As in the proof of Case~2 of the proof of Theorem~\ref{Theorem1.4}
this implies that equation~\eqref{e4.8} holds which in turn implies
$K$ and $K_0$ are homothetic.
\end{proof}

\providecommand{\bysame}{\leavevmode\hbox to3em{\hrulefill}\thinspace}
\providecommand{\MR}{\relax\ifhmode\unskip\space\fi MR }
\providecommand{\MRhref}[2]{%
  \href{http://www.ams.org/mathscinet-getitem?mr=#1}{#2}
}
\providecommand{\href}[2]{#2}

\end{document}